\documentclass{amsart}

\usepackage{amsmath}
\usepackage{amsthm}
\usepackage{amsfonts}
\title{Positivity conditions for Hermitian symmetric functions}
\author{John P. D'Angelo  and  Dror Varolin}

\newcommand{\noi}{\noindent}
\newcommand{\pf}{{\noi \it Proof. }}

\newcommand{\cg}{{\mathcal G}}
\newcommand{\ch}{{\mathcal H}}

\newcommand{\cp}{{\mathcal P}}

\newcommand{\C}{{\mathbb C}}

\begin{document}


\newtheorem{thm}{\sc Theorem}
\newtheorem{lem}{\sc Lemma}
\newtheorem{prop}{\sc Proposition}
\newtheorem{cor}{\sc Corollary}
\newtheorem{conj}{\sc Conjecture}
\newtheorem{defn}{\sc Definition}
\newtheorem{qn}{\sc Question}
\newtheorem{ex}{\sc Example}
\newtheorem{rmk}{\sc Remark}
\newtheorem{prop/defn}{\sc Proposition/Definition}

\noi {\sc Positivity conditions for Hermitian symmetric functions}\\
\medskip
{\it Dedicated to Yum-Tong Siu on the occasion of his sixtieth
birthday.}

\medskip

\noi {\sf John P. D'Angelo and Dror Varolin

\noi Dept. of Mathematics

\noi University of Illinois

\noi Urbana, IL 61801}

\section*{Introduction}

We introduce a family of positivity conditions for Hermitian
symmetric functions, establish basic properties, and connect the
ideas with complex geometry. Let $M$ be a complex manifold, and
let $M'$ denote its complex conjugate manifold. In this paper $M$
will typically be either $\C^n$ or the total space of a
holomorphic line bundle over a compact complex manifold.

A holomorphic function $R:M \times M' \to \C$ is called  {\it
Hermitian symmetric} if
$$ R(z, {\overline w}) = {\overline {R(w, {\overline z})}} $$
for all $z$ and $w$ in $M$. Observe that $z \to R(z, {\overline
z})$ is then necessarily a real-valued function; we say that $R$
is ``real on the diagonal''; conversely, by polarization $R$ is
determined by its values on the diagonal. Let $\cp_0(M)$ denote
the collection of Hermitian symmetric functions on $M \times M'$.
We will introduce in Definition 1, for $N$ a positive integer or
infinity, a subset $\cp_N = \cp_N(M)$ of $\cp_0(M)$. On the
diagonal, the set $\cp_1$ consists by definition of those $R$ that
are nonnegative, and the set $\cp_\infty$ consists (by a classical
result recalled in Lemma 2) of squared norms of Hilbert space
valued holomorphic mappings. The subsets $\cp_k$ therefore
interpolate two natural but distinct notions of nonnegativity for
Hermitian symmetric functions. The condition defining $\cp_k$ is
analogous to a positivity property for higher curvatures of
Hermitian metrics.

We show, as part of Theorem 3, that $\cp_j(\C^2) \ne \cp_k(\C^2)$
when $j \ne k$. Theorem 3 includes precise information about when
a member of a discrete collection of natural one-parameter
families of Hermitian symmetric polynomials lies in $\cp_k(\C^2)$.

Our other main result (Theorem 1) concerns a stability criterion.
A subset $S$ of $\cp_0(M)$ is called {\it stable} if there is a
finite $k$ for which $\cp_k \cap S = \cp_\infty \cap S$. The
minimum such $k$ is called the {\it stability index} of $S$, and
is written $I(S)$. For example, it is standard that the collection
$H$ of nonnegative Hermitian forms on $\C^n$  equals the
collection of squared norms of linear forms; therefore $I(H)=1$.
For each integer $N$ we give an example of an $S$ for which
$I(S)=N$. In Theorem 1 we relate the stability index to the number
of positive eigenvalues of the underlying matrix of coefficients
of a Hermitian symmetric function.

The concepts in this paper apply in a situation of some interest
in complex geometry. Let $X$ be a complex manifold and let $L$ be
a holomorphic line bundle over $X$. A (possibly degenerate) metric
$g$ on the fibres of $L$ is called a {\it globalizable singular
metric} if it is the restriction of a Hermitian symmetric function
$G$ on $L \times L'$. Our focus here will be on $G$. The subsets
$\cp_N$ of these metrics provide intermediate conditions between
simply being a (possibly degenerate) metric ($\cp_1$) and being a
holomorphic pullback of the Fubini-Study metric ($\cp_\infty$).

The proof of Theorem 1 applies to show that
the stability index for the collection $S$ of globalizable metrics on
a holomorphic line bundle $L$ over a compact complex manifold $X$
is bounded above by ${\rm dim} H^0(X,L^*) - 1$. Thus $S$ is stable. See Theorem 2.

There are many other useful positivity conditions for Hermitian
symmetric functions. In particular, the set of nonnegative
logarithmically plurisubharmonic (Hermitian) functions properly
contains $\cp_2$. We briefly discuss these issues in Section V.

The authors wish to acknowledge useful discussions over the years
with Yum-Tong Siu; in particular Siu pointed out to the first
author a connection between positivity conditions for
bihomogeneous polynomials and isometric imbedding. Motivated by
Siu's comments, the first author and Catlin \cite{cd2} proved an
isometric imbedding theorem for holomorphic vector bundles. This
theorem then motivated much of the work in the current paper. The
authors also wish to thank Dan Lichtblau of Wolfram Research, who
did a Groebner basis computation which provided useful insight,
and Dan Grayson for some useful discussions.

The first author acknowledges support from NSF grant DMS-0200551.

\section{Positivity Classes for Hermitian symmetric
functions}\label{I}

We begin by describing an essentially general class of Hermitian
symmetric functions on a complex manifold $M$. Let $\ch$ be a
Hilbert space with inner product $\langle \ , \ \rangle$. Let $f,g
:M \to \ch$ be holomorphic mappings, and define $R$ by the formula
\begin{eqnarray}\label{e1}
R(z, {\overline w}) = \langle f(z), f(w) \rangle - \langle g(z),
g(w) \rangle.
\end{eqnarray}
Then $R$ is Hermitian symmetric. When $g=0$ in (\ref{e1}) we have
$R=||f||^2$ on the diagonal; squared norms of holomorphic mappings
will play a special role in this paper. We also note that, by
choosing $\ch = \C$, and by choosing $f$ and $g$ appropriately in
(\ref{e1}), we can obtain the (polarization of the) real and
imaginary parts of an arbitrary holomorphic complex-valued
function.

When $R$ satisfies (\ref{e1}), we use the term {\it holomorphic
representation} to denote the holomorphic mapping $(f,g):M\to
\ch\times \ch$ determining $R$. When $g$ can be chosen to be $0$,
we simply say that $f:M \to \ch$ is a holomorphic mapping
representing $R$.

We now introduce the positivity classes $\cp_k(M)$ for the
collection of Hermitian symmetric functions on the manifold $M
\times M'$. Our notion evokes some classical functional analysis
such as Bochner's theorem on functions of {\it positive type}; the
key difference is that we focus on nonnegativity of matrices of a
fixed size rather than of all sizes.

\begin{defn}\label{d1}
(Positivity classes). Let $M$ be a complex manifold. We denote the
set of Hermitian symmetric function on $M \times M'$ by
$\cp_0(M)$.  For each positive integer $N$ we write $R \in
\cp_N(M)$ if
$$ \sum_{i,j=1}^N R(z_i, {\overline z_j}) a_i {\overline a_j} \ge 0 $$
for all $z= (z_1,...,z_N) \in M^N$ and all $a \in \C^N$. In other
words the Hermitian matrix with $(i,j)$ entry equal to $R(z_i,
{\overline z_j})$ is nonnegative definite. We write $R\in
{\cp}_\infty$ when $R\in \cp_N$ for all $N$; thus
$$ \cp_\infty = \bigcap_{N=1}^\infty \cp_N. $$
When $M$ is fixed, we drop it from the notation, and write $\cp_N$
for $\cp_N(M)$.
\end{defn}

\begin{rmk}\label{r1}
For each subset $\cp_N$ there are corresponding sharp forms; for
example we could demand that the matrix $R(z_i, {\overline z_j})$
be positive definite whenever the points are distinct. In some
contexts other sharp forms are useful; see Definition 5 in Section
V.
\end{rmk}

\begin{defn}\label{d2}
(Stability Index) Let $S$ be a subset of $\cp_0$. We define $I(S)$
to be the smallest $k$ for which
$$ S \cap \cp_\infty = S \cap \cp_k.$$
If no such $k$ exists we write $I(S) = \infty$. When $I(S)$ is finite we say that
$S$ is {\it stable}.
\end{defn}
We begin by noting some obvious properties of the sets $\cp_k$,
and then we study ${\cp}_\infty$.

\begin{lem}\label{l1}
Each $\cp_k$ is closed under sum and under product. For each $k$
we have $\cp_{k+1} \subset \cp_k$. If $R_\lambda$ is a family of
Hermitian symmetric functions depending continuously on some
parameter $\lambda$, then the set of $\lambda$ for which
$R_\lambda \in \cp_k$ is closed.
\end{lem}
\pf These facts follow easily from Definition 1. We note that the
proof of closure under product uses a well-known lemma of Schur:
if $(a_{ij})$ and $(b_{ij})$ are nonnegative definite matrices of
the same size, then their Schur product $( a_{ij} b_{ij})$ is also
nonnegative definite. \qed

\begin{lem}\label{l2}
Suppose that there is a Hilbert space $\ch$ and a holomorphic
function $f:M \to \ch$ such that
\begin{eqnarray}\label{e2}
R(z, {\overline w}) = \langle f(z), f(w) \rangle.
\end{eqnarray}
Then $R \in \cp_\infty$. Conversely, if $R \in \cp_\infty$, then
there is a Hilbert space $\ch$ and a holomorphic mapping $f:M \to
\ch$ representing $R$, and thus (\ref{e2}) holds.

\pf First assume (\ref{e2}) holds. Fix $N$, and choose arbitrary
points $z_i$ in $M$ and $a \in \C^N$. We see that
$$
\sum_{i,j=1}^N R(z_i, {\overline z_j}) a_i {\overline a_j} =
\sum_{i,j=1}^N \langle f(z_i)a_i , f(z_j) a_j \rangle = ||
\sum_{i=1}^N f(z_i)a_i||^2 \ge 0,
$$
and hence $R \in \cp_N$ for all $N$.
\end{lem}

The converse assertion is classical; various versions go back to
Mercer in 1909 and to E. H. Moore in 1916. For historical remarks
and related ideas we refer to \cite{a} and \cite{s}.

We give a sketch, following \cite{cw}, of the converse assertion.
Consider the complex vector space $V$ of functions on $M$ with
finite support. Using $R$ we define a Hermitian form on $V$ by the
following formula. For $u,v \in V$ we put
$$
\langle u,v\rangle_R = \sum R(z, {\overline w}) u(z) {\overline
{v(w)}},
$$
where the sum is taken over all $z, w \in M$, but is finite by the
support condition. Let $V_0$ denote the collection of $u$ with
$\langle u,u\rangle_R = 0$. It follows from the nonnegative
definiteness of all matrices $R(z_i, {\overline z_j})$ that $V_0$
is a subspace of $V$.  The quotient space is then an inner product
space under $\langle \ , \ \rangle_R$, and we take $\ch$ to be its
completion. For $z \in M$ we define $f(z)$ to be the image under
the quotient map of the function with support at the single point
$z$ and value unity there. Since $R$ is holomorphic on $M \times
M'$, it follows that $f: M \to \ch$ is holomorphic.\qed

\medskip

Next we give a simple example showing that $\cp_1(M)$ and
${\cp}_2(M)$ are distinct sets even in very simple situations.
This example also gives insight into the stability index.

\begin{ex}\label{ex1}
Let $M = \C^2$. For each real number $c$ we define $R_c$ by
$$
R_c(z,{\overline w}) = z_1^2 {\overline w_1}^2 + (c-2) z_1 z_2
{\overline w_1} {\overline w_2}  +z_2^2 {\overline w_2}^2
$$
It is elementary to check that $R_c \in \cp_1$ if and only if $c
\ge 0$ whereas $R_c \in \cp_2$ if and only if $c \ge 2$. In fact,
for each $k$ larger than $2$, $R_c \in \cp_k$ if and only if $c\ge
2$; therefore, if $S= \{ R_c :c \ge 0\}$, then $I(S)=2$. Recall by
contrast that $I(H)=1$ when $H$ is the set of nonnegative
Hermitian forms.
\end{ex}

We recall the standard test for nonnegative definiteness of an $N$
by $N$ matrix and provide the appropriate caution. A Hermitian
matrix is nonnegative definite if and only if every principal
minor determinant is nonnegative. It is not sufficient to assume
only that every {\it leading} principal minor determinant is
nonnegative. The matrix
$$
\left (
\begin{array}{ccc}
1 & 0 & 0 \\
0 & 0 & 0 \\
0 & 0 & -1
\end{array}
\right )
$$
has leading principal minor determinants of $1$, $0$, and $0$, and
yet it has a negative eigenvalue. See \cite{d4} for more
discussion.

We immediately obtain a method for deciding whether $R \in \cp_k$.
Given $k$ points $z_1, ..., z_k$  in $M$, we put $\Delta_k^R(z) =
{\rm det}(R(z_i, {\overline z}_j))$.

\begin{lem}\label{l3}
Suppose $R \in \cp_0(M)$ and $k \ge 1$. Then $R \in \cp_k$ if and
only if $R \in \cp_{k-1}$ and $\Delta_k^R(z) = {\rm det}(R(z_i,
{\overline z}_j)) \ge 0$ for all $z = (z_1,...,z_k) \in M^k$.
\end{lem}
Next we give a method for computing $\Delta_k^R(z)$ in many cases.
We say that $R$ has {\it finite rank} if there is a
finite-dimensional Hilbert space representation for $R$; in other
words, if there are holomorphic mappings $f,g : M \to \C^K$ such
that
\begin{eqnarray*}
R(z, {\overline w}) = \langle f(z), f(w) \rangle - \langle g(z),
g(w) \rangle.
\end{eqnarray*}
For example, every polynomial function $R$ has finite rank. See
\cite{d4} for a simple proof. One notes that globalizable metrics
(See Definition 3) on a holomorphic line bundle $L$ over a compact
complex manifold also have finite rank, simply because the space
of global sections of $L^*$ is finite-dimensional.

Suppose that $R$ has finite rank.
For each $z=(z_1,...,z_k) \in M^k$ we
consider a $k$ by $k$ matrix whose $i, j$ entry
is $h_j(z_i)$, where $h_j$ is a component function of either $f$ or $g$.
We let $\Delta_{k,m}^R (z)$ denote the tuple of determinants of all such matrices,
excluding the obvious repetitions,
for which precisely $m$ of the column vectors are components of $g$. We write
$||\Delta_{k,m} ^R (z)||^2$ for the sum of squared absolute values of all these
determinants.

We have the formula
\begin{eqnarray}\label{e3}
{\rm det}(R(z_i, {\overline z}_j)) = \Delta_k^R (z) = \sum_{m=0}^k
(-1)^m ||\Delta_{k,m}^R(z)||^2.
\end{eqnarray}

Formula (\ref{e3}) expresses the determinant $\Delta^R_k(z)$ in
terms of the components of $f$ and $g$ in a somewhat tractable
fashion. The proof is an elementary computation in exterior
algebra. One computes the wedge product of the column vectors of
the matrix ${\rm det}(R(z_i, {\overline z}_j))$, expands by the
distributive law, and uses the formula for the determinant of a
$k$ by $k$ matrix in terms of $k$-th exterior powers. Formula
(\ref{e3}) results. We note also that formula (29) from Section V
provides a proof and additional insight when $k=2$.

Formula (\ref{e3}) immediately combines with Lemma 3 to yield a
necessary and sufficient condition for an $R$ with finite rank to
lie in some $\cp_k$:

\begin{prop}\label{p1}
$R \in \cp_k$ if and only if $\sum_{m=0}^j (-1)^m
||\Delta_{j,m}^R(z)||^2 \ge 0$ for $1 \le j \le k$.
\end{prop}

We now turn to our first main result. Theorem 1 generalizes the
simple fact that, for scalar-valued entire holomorphic functions
$f$ and $g$ on $\C^n$, the function $|f|^2 - |g|^2$ can be
nonnegative only when $g$ is a constant multiple of $f$. For
general $k$, if $R = ||f||^2 - |g|^2 \in \cp_k$ and $f$ has too
few linearly independent components, then $g$ must be a linear
combination of the components of $f$, and thus $R$ actually must
be a squared norm.

Let $R$ be a Hermitian symmetric function on $\C^n$ with Taylor
expansion

$$ R(z,{\overline z}) = \sum_{a,b} c_{ab}z^a {\overline z}^b. $$
We let $N_+(R)$ and $N_-(R)$ denote the number of positive and negative eigenvalues
of the (perhaps infinite) Hermitian matrix $c_{ab}$.
We allow the value infinity.
In Theorem 1 we will not assume that $R$ has finite rank, although the situation
reduces easily to that case.

\begin{thm}\label{t1}
Let $R$ be a Hermitian symmetric function on $\C^n$ such that $R
\in \cp_k$. Then either $N_-(R) = 0$ or $N_+(R) \ge k+1$.
\end{thm}

\pf We will work solely on the diagonal. If $N_+(R) \ge k+1$, then
we are done. Suppose otherwise that $R\in \cp_k$, $N_-(R) > 0$,
and $N_+(R) \le k$. Since $N_+(R)$ is finite, we may write $R=
||f||^2 - ||g||^2$ as usual. In view of Lemma 1, by adding an
appropriate squared norm to $R$, we may suppose the following
hold: $g$ takes values in $\C$, $N_-(R)=1$, $f$ takes values in
$\C^k$, and that $N_+(R) = k$.

Choose $k$ points $z_1,...,z_k \in \C^n$. We define the square
matrix $A(f)(z)$ and column vector $G(z)$ by

$$  A(f)(z) = \left (
\begin{matrix}
f_1(z_1) &  f_2(z_1)  & ... & f_k(z_1) \cr f_1(z_2) & f_2(z_2) &
... & f_k(z_2) \cr ... & ... & & ... \cr f_1(z_k) & f_2(z_k) & ...
& f_k(z_k)
\end{matrix}
\right ) \ \ {\rm and} \ \ G(z) = \left (
\begin{matrix} g(z_1) \cr g(z_2) \cr ...\cr g(z_k)
\end{matrix}
\right ). $$

Consider the system of $k$ linear equations
$$A(f)(z) c(z) = G(z)$$
for $k$ unknowns $c_j(z)$. Since the functions $f_i$ are linearly
independent, the generic value of ${\rm det} (A(f)(z))$ is
nonzero. The solution to the system is therefore given as the
ratio of two determinants by Cramer's rule. We let $B_j(f,g)(z)$
denote the matrix obtained by replacing the $j$-th column in
$A(f)(z)$ with $G(z)$. Thus

\begin{eqnarray}\label{e4}
c_j(z) =\frac{{\rm det}(B_j(f,g)(z))}{{\rm det}(A(f)(z))}.
\end{eqnarray}
On the other hand, when $g$ is scalar-valued, (\ref{e3})
simplifies to give:
\begin{eqnarray}\label{e5}
{\rm det}(R(z_i, {\overline z}_j)) = \Delta_k^R (z) = |{\rm det}
(A(f)(z))|^2 - \sum_{j=1}^k |{\rm det}( B_j(f,g)(z))|^2
\end{eqnarray}
Since $R\in \cp_k$, the determinant on the left-hand side of (5)
is nonnegative. Each term in the sum on the right-hand side of (5)
is therefore bounded by $|{\rm det} A(f)(z)|^2$. Combining this
fact with (4) shows that each $c_j(z)$ is a bounded meromorphic
function, and hence a constant, written $c_j$. We therefore
obtain, for each $i$,
$$ g(z_i)  = \sum_j c_j f_j (z_i). $$
Since the points $z_i$ are arbitrary, $g$ is a linear combination
of the $f_j$; this statement contradicts the original assertion
that $N_-(R) \ne 0$. \qed

\begin{cor}\label{c1}
Let $S_k$ denote the collection of Hermitian symmetric functions
on $\C^n$ whose underlying matrices of Taylor coefficients have at
most $k$ positive eigenvalues. Then $I(S_k) \le k$ and $S_k$ is
stable.
\end{cor}

\medskip
The proof of Theorem 1 yields a test for
whether $g$ is a linear combination of the $f_j$. See Proposition 3
and Theorem 2 for an alternative point of view.

The proof of Theorem 1 goes through with essentially no change if
$\C^n$ is replaced with a holomorphic vector bundle over a compact
manifold, or, more generally, with a complex manifold admitting no
nonconstant bounded holomorphic functions. On the other hand, no
result resembling Theorem 1 holds for Hermitian symmetric
functions on bounded domains in $\C^n$. It is easy to write down,
for an arbitrary $j$, a Hermitian symmetric function $R$ on the
unit ball for which $N_+(R) = 1$ and $N_-(R) = j$.

\section{Globalizable metrics}\label{II}

Hermitian symmetric functions arise naturally in complex geometry.
In several important contexts Hermitian symmetric functions
(nonnegative on the diagonal) restrict to (possibly singular)
Hermitian metrics, and have played a key role in analytic
geometry. See \cite{siu} and \cite{cd2}.

Let $X$ be a compact, complex manifold, and suppose that $L$ is a
holomorphic line bundle over $X$. The
complex vector space $ H^0(X,L^*)$ of sections of the dual bundle $L^*$ is
finite-dimensional. Let $C$ be a Hermitian form on $ H^0(X,L^*)$.
We associate with $C$ the Hermitian symmetric function
$R_C$ on $L \times L'$ as follows.
If $\{\phi_\alpha \}$ is a basis for $H^0(X, L^*)$, and
$c_{\alpha \beta} = \langle C\phi_\alpha, \phi_\beta \rangle $, then
\begin{eqnarray}\label{e6}
 R_C(z, {\overline w}) =
\sum_{\alpha, \beta} c_{\alpha \beta} \phi_\alpha(z) {\overline
\phi_\beta(w)}.
\end{eqnarray}

\begin{defn}\label{d3}
Let $C$ be a Hermitian form on $ H^0(X,L^*)$. The function $R_C$
defined on $L \times L'$ by (6) is called a {\it globalizable
singular metric} on $L$.
\end{defn}

Suppose that $C$ is positive semi-definite; then, for each $\phi
\in H^0(X,L^*)$,  the function $R_C$ is nonnegative on the
diagonal. On the other hand, simple examples (see Example 1 and
Theorem 3) show that the function defined by (6) can be
nonnegative on the diagonal even when $C$ has some negative
eigenvalues. In case $C$ has negative eigenvalues the mapping $g$
in a holomorphic representation of $R_C$ must be nonzero. Hence
the intermediate positivity classes $\cp_k$ provide useful
geometric information.

A globalizable singular metric $G$
on $L$ can be written

\begin{eqnarray}\label{e7}
G = \sum_{\alpha, \beta} c_{\alpha \beta} \phi_\alpha {\overline
\phi_\beta},
\end{eqnarray}
where $\{\phi_\alpha \}$ form a basis for $H^0(X, L^*)$. The
(necessarily) Hermitian matrix $(c_{\alpha \beta})$ is called the
{\it underlying matrix} of $G$.

\medskip
We will drop the adjective {\it singular} and refer to $G$ as a
globalizable metric on $L$. We use the term {\it metric} because
$G$ restricts to a Hermitian metric on the fibres of $L$. In other
words, there is a (possibly singular) metric $g$ in the usual
sense for which $G(z, {\overline w}) = g(z,{\overline w})$
whenever $\pi(z) = \pi(w)$. In \cite{cd2} the metric $g$ is
defined to be globalizable when a $G$ extending $g$ and satisfying
(7) exists.

Complex projective space ${\bf P}_{n-1}$ and powers of the universal bundle
${\bf U}^m$ provide nice examples.
The usual Fubini-Study metric $G$ on ${\bf U}$ is globalizable.

The dual bundle ${\bf H}^m$ of ${\bf U}^m$ is the $m$-th power of
the hyperplane bundle; it is generated by global sections which we
may identify with homogeneous polynomials of degree $m$. It is
natural to equip ${\bf U}^m$ with the Hermitian metric $G^m$ given
by the $m$-th tensor power of the Fubini-Study metric on ${\bf
U}$. The metric $G^m$ provides another simple example of a
globalizable metric. In the particular case of ${\bf U}^m$, a
globalizable metric may be identified with a bihomogeneous
polynomial on $\C^n$:

$$ R(z,{\overline w}) = \sum_{|\alpha|=m, |\beta|=m}
c_{\alpha\beta} z^\alpha {\overline w}^\beta. $$
The matrix of coefficients $(c_{\alpha \beta})$ is Hermitian symmetric
and $R(z, {\overline z})$ nonnegative.
For ${\bf U}^m$,
the natural metric $G^m$ can be identified with the bihomogeneous
polynomial defined by
$$ R(z,{\overline w}) = \langle z, w \rangle^m. $$
In Theorem 3 we will consider perturbations of this particular $R$.

\medskip

We interpret one piece of Lemma 1 in this setting.

\begin{cor}\label{c2}
Suppose that $(L,R)$ and $(E,G)$ are holomorphic line bundles over
a complex manifold $X$ with the indicated globalizable metrics $R$
and $G$. Suppose that $R$ and $G$ are in $\cp_k$. Then $R \otimes
G \in \cp_k$.
\end{cor}
\pf The result follows from Lemma 1, because $\cp_k$ is closed
under product. \qed

\medskip

Corollary 2 implies that the natural tensor product metric on a
power of a line bundle lies in the same class $\cp_k$ as does the
original metric. When $k=2$ the converse of this fact plays a
crucial role in the proof of the isometric imbedding theorem in
\cite{cd2}.

\section{Bihomogeneous polynomials, metrics, and the functions
$\Sigma^k$}\label{III}

In this section we give an alternative manner for verifying that
$R \in \cp_k$ when $R$ is a bihomogeneous polynomial, or more
generally, a globalizable metric. We begin with a general lemma.

\begin{lem}\label{l4}
Let $M$ be a complex manifold and $\ch$ be a Hilbert space with
inner product $\langle \ , \ \rangle$. Suppose that $f,g :M \to
\ch \times \ch$ is a Hilbert space representation of the Hermitian
symmetric function $R$. Then $R \in \cp_N$ if and only if
\begin{eqnarray}\label{e8}
|| \sum_{i=1}^N f(z_i)a_i||^2 \ge || \sum_{i=1}^N g(z_i)a_i||^2
\end{eqnarray}
for all choices of $N$ points $z_i \in M$ and all $a \in \C^N$.
\end{lem}

\pf The computation is virtually the same as the proof of Lemma 2,
and hence is left to the reader. \qed
\medskip

We recall that a bihomogeneous polynomial $R$ on $\C^n \times
(\C^n)'$ is a polynomial in $z$ and ${\overline w}$ that is
homogeneous of the same degree in both sets of variables.
Equivalently, for $\lambda \in \C$,
\begin{eqnarray}\label{e9}
R(\lambda z, {\overline {\lambda z}}) = |\lambda|^{2m} R(z,
{\overline z}).
\end{eqnarray}

A bihomogeneous polynomial is real on the diagonal if and only if
it is Hermitian symmetric, and this symmetry holds if and only if
the matrix of coefficients of $R$ is Hermitian symmetric. See
\cite{d2}and \cite{d4} for discussion of bihomogeneous
polynomials.

We may identify a Hermitian symmetric nonnegative bihomogeneous
polynomial $R$ of degree $2m$ on $\C^n$ with a globalizable metric
on ${\bf U}^m$ over ${\bf P}_{n-1}$. To see this fact, we write
$R$, as in the definition of globalizable metric, in the form (7),
where $\phi_\alpha$ is the monomial $z^\alpha$. The dual bundle
${\bf H}^m$ is generated by global sections which we may identify
with homogeneous polynomials of degree $m$.

Let now $X$ be a compact complex manifold and $\pi: L \to X$ a
holomorphic line bundle. We denote the pairing of a section $s$ of
$L^*$ with a vector $v \in L$ by $s[v]$. We obtain Hermitian
symmetric functions $R: L \times L' \to \C$ by mimicking the
situation for bihomogeneous polynomials.

Every Hermitian symmetric function $R$ on $L \times L'$ is given,
after diagonalizing the associated Hermitian form $C$ from (6) by
the formula
\begin{eqnarray}\label{e10}
R(v,{\overline w}) = \sum_{j=1}^m \mu_j s_j(\pi(v)) [v] \
{\overline {s_j(\pi(w)) [w] }}.
\end{eqnarray}
In (10) the $\mu_j$ are the nonzero eigenvalues of $C$ and
$s_1,...,s_k$ are linearly independent elements of $H^0(X, L^*)$.
By collecting the terms according to the sign of $\mu_j$ we write
$R = ||f||^2 - ||g||^2$, where the components of $f$ and $g$ are
global sections of $L^*$. Thus the sections determine a
(finite-dimensional) Hilbert space representation for $R$, and $R
\in \cp_1$ precisely when $R$ is a globalizable metric. The
following simple result characterizes when $R \in \cp_k$.

\begin{prop}\label{p2}
Let $\pi: L \to X$ be a holomorphic line bundle over a compact
complex manifold $X$. Let $R$ be the Hermitian symmetric function
defined on the diagonal by $R=||f||^2 - ||g||^2$. Then $R \in
\cp_k$ if and only if
\begin{eqnarray}\label{e11}
|| \sum_{i=1}^k f(w_i)||^2 \ge || \sum_{i=1}^k g(w_i)||^2
\end{eqnarray}
for all $w_1, ..., w_k$ in $X$.
\end{prop}
\pf Let $s$ be a section of $L^*$, and choose $v \in L$. For $a
\in \C$ and $p \in X$, we have $ s(p) [a v]= s(p) [w]$ for some $w
\in L$. The proposition follows by combining this fact with
formula (8) and the conclusion of Lemma 2. \qed

We introduce notation for the functions appearing
inside the norms in (11).
\begin{defn}\label{d4}
Let $X$ be a set, $G$ be an additive Abelian group, and $f:X\to G$
a function. For each positive integer $N$ we define $\Sigma^N f :
X^N \to G$ by
$$
(\Sigma^N f)(z_1, ..., z_N) = \Sigma_{i=1}^N f(z_i).
$$
\end{defn}
We can then interpret the conclusion of Proposition 2 in several
ways.

\begin{cor}\label{c3}
Let $R$ be a Hermitian symmetric function on $L\times L'$, written
in the form $R= ||f||^2 - ||g||^2$ as in Proposition 2. Then $R
\in \cp_N$ if and only if $|| \Sigma^N f||^2 \ge || \Sigma^N
g||^2$.
\end{cor}

\begin{cor}\label{c4}
Let $R_\lambda$ be defined on the diagonal of $L \times L'$ by
$R_\lambda=||f||^2 - \lambda ||g||^2$. Let $S= \{R_\lambda\}$. If
there is an $N$ such that the variety $V(\Sigma^N f)$ is not
contained in the variety $V(\Sigma^N g)$, then $I(S) \le N$.
\end{cor}

Let $X$ be a compact complex manifold, and let $L \to X$ be a
holomorphic line bundle.  We let $h = {\rm dim} \ H^0(X,L^*)$, By
Definition 4, $\Sigma ^k s : L^k \to \C$ is defined by
$$(\Sigma ^k s) (v_1,...,v_k) = \sum _{j=1} ^k s[v_j],$$
and is a section of the vector bundle $L^k \to X^k$.  We will consider
the zero varieties
$$V(\Sigma ^k s) := \left \{ (x_1,...,x_k) \in X^k \ ;\ \Sigma ^k s
|L_{x_1} \times ... \times L_{x_k} \equiv 0 \right \}.$$

\begin{prop}\label{p3}
Let $g,s_1,...,s_N \in H^0 (X,L^*)$, and suppose that
$$\bigcap _{j=1} ^N V(\Sigma ^h s_j) \subset V(\Sigma ^h g).$$
Then $g$ is a linear combination of $s_1,...,s_N$.
\end{prop}

\begin{lem}\label{l5}
Let $T \in H^0(X,L^*)^*$.  Then there exist points $v_1,...,v_h
\in L$ such that for all $s \in H^0(X,L^*)$,
$$ Ts = \Sigma ^h s(v_1,...,v_h).$$
\end{lem}

\pf We may assume $T \neq 0$.  Choose a basis $s^1,...,s^h$ of
$H^0(X,L^*)$ such that
$$Ts^1 = 1 \quad {\rm and} \quad Ts^j = 0 \ {\rm for}\ j=2,...,h.$$
Select vectors $u_1,...,u_h$ such that the matrix with components
$A^i_j := s^i(u_j)$ is invertible.  This choice
is trivial (as is the Lemma) when $h=1$.
For $h \ge 2$ selecting such vectors is also possible; for instance, we may choose
$$u_j \in V(s^j) - \bigcup _{i \neq j} V(s^i), \quad j=1,...,h.$$
Now choose complex numbers $\lambda _1,...,\lambda _h$ such that
$$\sum _{j=1} ^h A^1_j\lambda _j = 1 \quad {\rm and} \quad \sum _{j=1} ^h
A^i_j\lambda _j = 0 \ {\rm for}\ i=2,...,h,$$ and let $v_j =
\lambda _j u_j$, $j=1,...,h.$  For any $s \in H^0(X,L^*)$, there
exist $\mu _i$ such that $s = \sum _i \mu _i s^i$.  In this case,
$Ts = \mu _1$.  On the other hand,
\begin{eqnarray}\label{e12}
\qquad \Sigma ^h s (v_1,...,v_h) = \sum _{i=1} ^h \sum _{j=1} ^h
\mu _i s^i(v_j) = \sum _{i=1} ^h \sum _{j=1} ^h \mu _i A^i _j
\lambda _j = \sum _{i=1} ^h  \mu _i \delta ^{1i} = \mu _1.
\end{eqnarray}
Formula (12) completes the proof of the Lemma. \qed

\medskip

\noi {\it Proof of Proposition 3.}  Let $U$ be the subspace of
$H^0(X,L^*)$ generated by the sections $s_1,...,s_N$.  Let
$U^{\perp} \subset H^0(X,L^*) ^*$ be the annihilator of $U$.
Suppose $T \in U^{\perp}$. Then by Lemma 5 there exist vectors
$v_1,...,v_h \in L$ such that $Ts = \Sigma ^h s(v_1,...,v_h)$.
Now, since $Ts_i = 0$ for $i=1,..., N$, the hypothesis on $g$ and
the expression $Ts = \Sigma ^h s(v_1,...,v_h)$ imply that $Tg =0$.
Thus $g \in (U^{\perp}) ^{\perp} = U$, as desired. \qed

\medskip
It is possible to prove Proposition 3 using the ideas in Theorem 1, but
we feel that Proposition 3 is interesting in its own
right.  The operation $\Sigma ^k$ can be thought of as
a discrete version of
the $(k-1)$--jet of a section.  In fact, the $\Sigma ^k$ contain more
information than jets do, because one may look at
points not on the diagonal. From this point of view the hypothesis of Proposition 3
says that the global section $g$ vanishes to order $h$ whenever
the global sections $s_1,...,s_N$ do.

\medskip
In both Theorem 1 and Proposition 3 we conclude that a function
$g$ must be a linear combination of some given functions $f_j$ (or
$s_j$) under certain conditions. In each case these conditions
guarantee that the function $ ||f||^2- |g|^2 $ lies in
$\cp_\infty$. Thus Proposition 3 is an algebraic analogue of Lemma
2.

Let $S$ denote the set of singular globalizable metrics on $L$.
Proposition 3 and Corollary 3 combine to show that
$I(S) \le {\rm dim} H^0(X,L^*)$.
The ideas from Theorem 1 provide a stronger
result on the stability index.

\begin{thm}\label{t2}
(Stability of Positivity Conditions) Let $X$ be a compact complex
manifold, and let $L$ be a holomorphic line bundle over $X$. Let
$S$ be the set of globalizable metrics on $L$. Then $S$ is stable.
In fact $I(S)$ is at most ${\rm dim}( H^0(X,L^*)) -1$.
\end{thm}
\pf  We first observe that bounded meromorphic functions on $L$
must be constant. We can therefore virtually repeat the proof of
Theorem 1. Put $k= {\rm dim}( H^0(X,L^*)) -1$. Let $R$ denote a
globalizable metric, whose underlying Hermitian form is $C$. Note
that $C$ has at most $k+1$ eigenvalues. Suppose that $R \in
\cp_k$; by the argument in Theorem 1, either $C$ has no negative
eigenvalues, or at least $k+1$  positive eigenvalues. In the
second case the total number of positive eigenvalues must equal
$k+1$. Thus, in either case, $R$ is in $\cp_\infty$. Therefore $S
\cap \cp_k = S \cap \cp_\infty$, and hence $S$ is stable and
$I(S)=k$. \qed

\section{Perturbations of the standard metric on ${\bf U}^{2m}$.}
\label{IV}

In this section we will show that the sets $\cp_k(M)$ are distinct
even when $M = \C^2$. We may also regard formula (13) below as
arising from a family of globalizable metrics on the line bundle
${\bf U}^{2m}$ over ${\bf P}_1$.

We consider a family of bihomogeneous polynomials defined on
$\C^2$. These polynomials may be considered as perturbations of
the natural example (on the diagonal) $z \to ||z||^{4m} = R_0(z,
{\overline z})$. For each real number $\lambda$, and each positive
integer $m$ we define
\begin{eqnarray}\label{e13}
R_\lambda (z,{\overline w}) = \langle z,w\rangle^{2m} - \lambda
(z_1z_2 {\overline w_1} {\overline w_2})^m
\end{eqnarray}
Note that $R_\lambda$ depends on both $\lambda$ and $m$. Let us
write $S_{\lambda,m}$ for the set of polynomials of the form (13).
Theorem 1 below provides precise information on when $R_\lambda
\in \cp_k$. We determine the set of values of $\lambda$ (in terms
of $m$) for this inclusion to hold in the special cases $k = 1$,
$k=2$, and $k=\infty$. We also show that $\cp_k = \cp_j$ for $k
\ge j \ge m+1$. Thus, $R_\lambda \in \cp_\infty$ if and only if
$R_\lambda \in \cp_k$ whenever $k \ge m+1$. In particular
$I(S_{\lambda,m}) = m+1$. Finally we determine the precise value
of $\lambda$ for $R_\lambda$ to be in $\cp_m$.

These results demonstrate an important fact
about our positivity conditions.

\medskip

\noindent {\sc Consequence}. {\it For each $m$, we have
${\cp}_m(\C^2 \ne\cp_{m+1}(\C^2)$. Thus $\cp_1(M)$ is not stable
already for $M = \C^2$.}

\begin{thm}\label{t3}
Fix $m$, and let $R_\lambda$ be defined by (13). The following
statements hold:
\begin{enumerate}
\item[1)] $R_\lambda \in \cp_1$ if and only if $\lambda \le 2^{2m}$.

\item[2)]  $R_\lambda \in \cp_2$ if and only if $\lambda \le
2^{2m-1}$.

\item[3)] $R_\lambda \in \cp_\infty$ if and only if
$\lambda \le \left ( ^{2m} _{\ m} \right ).$

\item[4)] For each integer $k >m$, $R_\lambda \in \cp_k$ if and only
if $R_\lambda \in \cp_\infty$.

\item[5)] $R_\lambda \in \cp_m$ only if $\lambda \le \left ( ^{2m} _{\ m}
\right ) +2$.

\item[6)] $I(S_{\lambda,m})=m+1$.
\end{enumerate}
\end{thm}

\noi {\it Proofs.}

\noi 1) To decide whether $\cp_1$ holds is elementary. By
homogeneity it suffices to assume that $||z||^2=1$; we then want
the largest $\lambda$ such that
\begin{eqnarray}\label{e14}
1 - \lambda |z_1|^{2m} |z_2|^{2m} \ge 0
\end{eqnarray}
given that $|z_1|^2 + |z_2|^2 = 1$. It is evident that the maximum
value of $|z_1|^{2m} |z_2|^{2m}$ on the sphere occurs when
$|z_1|^2 = |z_2|^2 = \frac{1}{2}$. Combining this observation with
(14) shows that
$$1 - \frac{\lambda}{2^{2m}} \ge 0,$$
which gives the result we want.

\medskip

\noi 3) To decide when $\cp_\infty$ holds is also elementary. On
the diagonal we have
$$
R_\lambda(z, {\overline z})=(|z_1|^2+|z_2|^2)^{2m}-\lambda
(|z_1|^2 |z_2|^2)^m.
$$
Expanding the left-hand side by the binomial theorem reveals that
we obtain a squared norm if and only $\lambda$ is at most the
coefficient there of $(|z_1|^2 |z_2|^2)^m$; this coefficient is
$\left ( ^{2m} _{\ m}\right )$.

\noi 2) We next verify that $\cp_2$ fails for $R_\lambda$ when
$\lambda > 2^{2m-1}$. We choose the points $z=(1,1)$ and
$w=(1,-1)$ and compute
\begin{eqnarray}\label{e15}
R_\lambda  (z, {\overline z}) R_\lambda  (w, {\overline w}) -
|R_\lambda  (z, {\overline w})|^2
\end{eqnarray}
there. We suppose that the expression in (15) is nonnegative.
Using the orthogonality of $z$ and $w$, the computation simplifies
and we obtain the condition
$$ 0 \le 2^{2m} - 2 \lambda, $$
which shows that the inequality fails for $\lambda > 2^{2m-1}$. To
verify that $\cp_2$ holds for $\lambda \le 2^{2m-1}$ one must
first show that the cut-off value of $\lambda$ is determined when
$\langle z, w\rangle = 0$. We omit the considerable details; they
are similar to the proof of 5). After verifying this fact, we must
maximize $|z_1 z_2|^{2m} + |w_1 w_2|^{2m}$ given that $|z_1|^2 +
|z_2|^2 = |w_1|^2 + |w_2|^2 = 1$ and that $z$ and $w$ are
orthogonal. The maximum occurs at several points, including
$\frac{1}{\sqrt{2}} (1,1)$ and $\frac{1}{\sqrt{2}} (1,-1)$. These
points yield the desired statement about $\lambda$.

We prove 4) in Propositions 4 below. We prove 5) in Proposition 5
below. Statement 6) follows by combining statements 4) and 5).
\qed

\medskip
Before proving these statements we interpret them when $m=1$ and $m=2$.
First we set $m=1$ and recover information implied by Theorems 1 and 2:
$$ \cp_1 \cap S_{\lambda,1} \ne \cp_2 \cap S_{\lambda,1} $$
but that
$$ \cp_2 \cap S_{\lambda,1} = \cp_\infty \cap S_{\lambda,1}. $$
Thus $I(S_{\lambda,1}) = 2$.

Next consider $m=2$; we see that
$$ \cp_1 \cap S_{\lambda,2} \ne \cp_2 \cap S_{\lambda,2}
 \ne \cp_3 \cap S_{\lambda,2} $$
but that
$$ \cp_3 \cap S_{\lambda,2} = \cp_\infty \cap S_{\lambda,2}. $$
Thus $I(S_{\lambda,2}) = 3$. This information does not follow from
Theorems 1 and 2.

It is instructive to prove 4) of Theorem 3 first when $m=2$. By
Corollary 3 we must determine the set of $\lambda$ for which the
inequality
\begin{eqnarray}\label{16}
&& |z_1 ^4 + w_1^4 + u_1^4|^2 +4 |z_1^3 z_2 + w_1^3 w_2 + u_1^3
u_2|^2 \\
\nonumber && \qquad \qquad + 4 |z_1 z_2^3 + w_1 w_2^3 + u_1 u_2^3
|^2 +|z_2 ^4 + w_2^4
+ u_2^4|^2\\
\nonumber && \qquad \qquad \qquad \qquad \ge (\lambda -6)
|(z_1z_2)^2 +(w_1w_2)^2 + (u_1u_2)^2|^2
\end{eqnarray}
holds for all triples of points $z,w,u$.

(16) is trivial for $\lambda \le 6$;
of course $6$ is the cut-off point for being a squared norm.
We use Corollary 4. Suppose that (16) holds,
and we can show that the left-hand side of (16) can
vanish when the right-hand side does not vanish. We see that the inequality
will fail for any $\lambda$ larger than $6$.

Let $h$ denote the polynomial in six variables given by
$(z_1z_2)^2 +(w_1w_2)^2 + (u_1u_2)^2$. Let $J$ denote the ideal
in the polynomial ring given by
$$ J = (z_1 ^4 + w_1^4 + u_1^4, z_1^3 z_2 + w_1^3 w_2 + u_1^3 u_2,
z_1 z_2^3 + w_1 w_2^3 + u_1 u_2^3, z_2 ^4 + w_2^4 + u_2^4). $$
To verify the existence of points as in the previous paragraph, it suffices
to prove that $h$ is not in the radical of $J$. Dan Lichtblau of Wolfram Research
verified this and other statements
for us using the Groebner basis algorithm from Mathematica.
After knowing this information, it is natural to instead seek points that work.

First we choose $z_1 z_2 = w_1w_2 = u_1u_2 = 1$.
With these choices $h$ will have the value $3$.
Next we assume that $|z_1|^2 = |w_1|^2 = |u_1|^2 = 1$.
By doing so, the four equations defined by the vanishing of the generators
of $J$ become the two equations (17) and (18) and their complex conjugates.
\begin{eqnarray}\label{e17}
z_1^4 + w_1^4 + u_1^4 =  0 \\
\label{e18} z_1^2 + w_1^2 + u_1^2 = 0
\end{eqnarray}
We choose the three values to equal $1$, $b$, and $b^2$, where $b$
is a primitive third root of unity. Then both (17) and (18) become
$1+b +b^2 = 0$, which holds since $b^3=1$ and $b\ne 1$. This
verifies statement 4) of Theorem 3 when $m=2$.

We now generalize this proof.
\begin{prop}\label{p4}
(Stability) For each positive integer $m$, and each real
$\lambda$, let $S_{\lambda,m}$ denote the set of Hermitian
symmetric functions $R_\lambda$ on $\C^2$ satisfying (19):
\begin{eqnarray}\label{e19}
R_\lambda (z,{\overline w}) = \langle z,w\rangle^{2m} - \lambda
(z_1z_2 {\overline w_1} {\overline w_2})^m.
\end{eqnarray}
Then, for $k \ge m+1$, the sets $S_{m,\lambda} \cap \cp_k$ are all
the same. In particular, $R_\lambda \in \cp_k$ if and only if
$\lambda \le \left ( ^{2m} _{\ m} \right )$, and $I(S_{\lambda,m})
\le m+1$.
\end{prop}
\pf We mimic the proof in the special case $m=2$. Let $g(z) =
z_1^m z_2^m$. We observe that, for $z_1 \ne 0$, we have
$g(z_1,\frac{1}{z_1}) = 1$. For any collection of points $W_j$ in
$\C^2$ of the form $(w_j, \frac{1}{w_j})$ we therefore have
$\Sigma^N g(W_1,...,W_N) = N$.

Let $J$ denote the ideal analogous to the ideal in the special
case above. The generators of $J$ are the functions $\Sigma^N
h_j$, where $h_j(z) = z_1^{2m-j}z_2^j$, for $ 0 \le j < m$ or $m+1
\le j \le 2m$. We claim that, for $N = m+1$, we can find $N$
points $W_i$ in $\C^2$ of the form $W_i = (w_i, \frac{1}{w_i})$
such that $\Sigma^N f(W_1,...,W_N) = 0$. By Corollary 4 this
implies the desired result for $R_\lambda$.

We will find solutions where
$|w_i|^2 =1$. Assume this condition. Then the $2m$ equations
$h_j=0$ become $m$ equations and their conjugates.
It therefore suffices to verify the first $m$ equations:
\begin{eqnarray}\label{e20}
\sum_{i=1}^m w_i^{2m-2j} = 0
\end{eqnarray}
for $0 \le j \le m-1$. To satisfy the equations in (20) we set
$(w_1)^2=1$, $(w_2)^2= \eta$, $(w_3)^2 = \eta^2$, and so on, where
$\eta$ is a primitive $m+1$-st root of unity. We obtain the $m$
equations
\begin{eqnarray}\label{e21}
\sum_{i=1}^m \eta^{(i-1)(m-j)} = \sum_{i=0}^{m-1}(\eta^{m-j})^i=0.
\end{eqnarray}
For each $j$, equation (21) holds because $\eta^{m-j}$ is an
$m+1$-st root of unity but does not equal $1$. \qed

\begin{prop}\label{p5}
Statement 5) from Theorem 3 holds.
\end{prop}

\pf We have $R_\lambda = ||f||^2 - \lambda||g||^2$ on the
diagonal, where $f$ and $g$ are as follows:
\begin{eqnarray}\label{e22}
||f(z)||^2 = ||z||^{4m} - \left ( ^{2m}_{\ m}\right ) |(z_1z_2)^m|^2\\
\label{e23} ||g(z)||^2 =  |(z_1 z_2)^m|^2.
\end{eqnarray}
Note that the right-hand side of (22) is a squared norm. The
components of $f$ are all the homogeneous polynomials of degree
$2m$ except for the middle term $(z_1z_2)^m$, and the coefficients
are the binomial coefficients.

By Corollary 3 the largest $\lambda$ for which $\cp_k$ holds is
given by the infimum of the expression
\begin{eqnarray}\label{e24}
\frac{||\Sigma^k f(w_1,...,w_k)||^2}
{||\Sigma^kg(w_1,...,w_k)||^2}.
\end{eqnarray}
In (24), $w_j = ( (w_j)_1, (w_j)_2 )$. Let $k=m$. Since (24) is
homogeneous of degree $0$, and because of the symmetry, the
infimum happens where $w_k = (1, \eta^{k-1})$ and $\eta^m=1$.
Plugging this into (24), and using properties of roots of unity,
all the terms in $||\Sigma^m f||^2$ vanish except the terms
$|\Sigma ((w_j)_1)^m|^2$ and $|\Sigma ((w_j)_2)^m|^2$. These each
give the value $m^2$. On the other hand, the value of $\Sigma^m g$
at these points is $m$. Hence the value of (24) at these points is
$2$. Incorporating the terms $\left ( ^{2m}_{\ m}\right )
|(z_1z_2)^m|^2$ gives the desired result for $R_\lambda$, where
$f$ and $g$ are expressed in the forms (22) and (23).

This calculation proves that $R_\lambda$ is not in $\cp_m$ when
$\lambda > \left ( ^{2m}_{\ m}\right ) + 2$. Containment does hold
for this value, because the infimum in (24) is attained at these
points. Even though the numerator vanishes on a two-dimensional
variety, the denominator vanishes there as well, and the limiting
value of the ratio, as we approach this variety, exceeds $2$. \qed

\begin{rmk}\label{r2}
For any homogeneous polynomial mapping $f$ of degree $m$ in two
variables, the variety $V(\Sigma^k f)$ is positive dimensional for
$k \ge 2$. On the other hand, for the functions in Theorem 3, we
have $V(\Sigma^k f) \subset V(\Sigma^k g)$ for $k \ge 2$. The key
point is that the infimum of the ratio happens away from the
variety $V(\Sigma^k f)$.
\end{rmk}

\section{The role of $\cp_2$}\label{V}

The main purpose of this section is to discuss $\cp_2$ in more
detail. For completeness we first recall and augment our list of
geometric positivity conditions. Suppose that $R=R_C$ is a
globalizable metric as in (6). The following positivity conditions
all arise:

\begin{defn}\label{d5}
(Geometric positivity conditions)
\begin{enumerate}
\item[1)] ${\cg \cp}_\infty$: The matrix $(c_{\alpha\beta})$ of
coefficients is nonnegative definite.

\item[2)] ${\cg \cp}_\infty^\sharp$: The matrix $(c_{\alpha\beta})$
of coefficients is positive definite.

\item[3)] ${\cg}$: There is an integer $d$ such that
$$R^d(z,w)=\sum_{\mu, \nu} E_{\mu \nu} \psi_\mu {\overline \psi}_\nu$$
where the matrix $(E_{\mu \nu})$ is nonnegative definite.

\item[4)] ${\cg}^\sharp$: There is an integer $d$ such that
$$R^d(z,w)=\sum_{\mu, \nu}E_{\mu \nu}\psi_\mu {\overline \psi}_\nu$$
where the matrix $(E_{\mu \nu})$ is positive definite.

\item[5)] $R \in \cp_2$; equivalently the global Cauchy-Schwarz
inequality holds:
$$|R(z,{\overline w})|^2 \le R(z,{\overline z}) \ R(w,{\overline w}).$$

\item[6)] A sharp form of the global Cauchy-Schwarz inequality holds:
$$ |R(z,{\overline w})|^2 \le R(z,{\overline z}) \ R(w,{\overline w}),$$
and equality happens only in some specified precise setting.

\item[7)] The function $z \to {\rm log} R(z, {\overline z})$ is
plurisubharmonic.

\item[8)] The function $z \to R(z, {\overline z})$ is plurisubharmonic.

\end{enumerate}
\end{defn}

\medskip
Items 2), 4) and 6) are sharp forms of items 1), 3), and 5). We could of course
also introduce sharp forms of 7) and 8).

By elementary linear algebra 1) is equivalent to the existence of
a (finite dimensional) Hilbert space valued mapping $f$
representing $R$ and is thus is the same as $\cp_\infty$. One can
easily show \cite{d4} that 3) implies 5). By definition, 5) is the
same as $R \in \cp_2$. The assumption that $R$ is a (possibly
degenerate) metric is of course the same as $R \in \cp_1$.
Conditions 7) and 8) lie between $\cp_1$ and $\cp_2$; It is shown
in \cite{d4} that 5) implies 7) implies 8). We revisit our main
example below to show that each converse assertion fails.
Statements 7) and 8) are equivalent in the bihomogeneous case; see
\cite{d4}.

\medskip

\noi {\sc Theorem 3 Revisited}. {\it Let $M$ be $\C^2$. For each
real number $\lambda$ we define a Hermitian symmetric polynomial
$r_\lambda$ by
\begin{eqnarray}\label{e25}
r_\lambda (z, {\overline z}) = (|z_1|^2 + |z_2|^2)^4  - \lambda
|z_1z_2|^4.
\end{eqnarray}
(Thus $m=2$ in Theorem 3.) The following statements are true:

\begin{enumerate}
\item[2.1)] $r_\lambda \in \cp_1$ if and only if $\lambda \le 16$. (This
is the condition for being a metric.)

\item[2.2)] $r_\lambda$ is plurisubharmonic (on the diagonal) for
$\lambda \le 12$.

\item[2.3)] ${\rm log}(r_\lambda)$ is plurisubharmonic (on the
diagonal) for $\lambda \le 12$.

\item[2.4)] $r_\lambda \in \cp_2$ if and only $\lambda\le 8$.

\item[2.5)] For $k \ge 3$, $r_\lambda \in \cp_k$ if and only if
$r_\lambda \in \cp_\infty$; this condition occurs if and only if
$\lambda \le 6$.
\end{enumerate}
}

\medskip
{\bf Remark 3}. Theorem 3 therefore reveals that four of the five
conditions are distinct already for bihomogeneous polynomials of
degree $8$ in $\C^2$. If we set $z_2=1$, and consider the
resulting function of one variable, then the condition $\lambda
\le 12$ for logarithmic plurisubharmonicity is unchanged; the
condition for plurisubharmonicity becomes $\lambda  \le
\frac{3}{32}(69 + 11\sqrt{33})$. (See \cite{d4}). The value of
this expression is approximately 12.39. Thus all five conditions
are distinct for polynomials in two variables.
\medskip

Condition 7) from Definition 5 is equivalent to the strict
negativity of the bundle $L$ using the metric $R$. The
Cauchy-Schwarz inequality on a globalizable metric is thus an
intermediate condition between  being a holomorphic pullback and
having negative curvature. We therefore give a different
characterization of $\cp_2$; this result holds for Hermitian
symmetric functions not necessarily arising from metrics.

In order to state the result in a geometric fashion we introduce
two pieces of notation. Suppose that $f,g:M \to \ch$ are
holomorphic mappings. We write $f \wedge g$ for the skew-symmetric
function from $M \times M$ to $\ch \otimes \ch$ defined by
\begin{eqnarray}\label{e26}
(f \wedge g) (z,w) = f(z) \otimes g(w) - f(w) \otimes g(z)
\end{eqnarray}
We also write $\theta_f(z,w)$ for the angle between the vectors
$f(z)$ and $f(w)$ in $\ch$. Finally we note the standard identity
\begin{eqnarray}\label{e27}
&& ||f(z) \otimes f(w) ||^2 = ||f(z)||^2 ||f(w)||^2 \\
\nonumber && \qquad \qquad  = ||f(z)||^2 ||f(w)||^2 {\rm sin}^2
\theta_f(z,w) + |\langle f(z),f(w)\rangle|^2.
\end{eqnarray}

Let $(f,g):M \to \ch \times \ch$ be a holomorphic representation
of $R$. We have the following beautiful geometric interpretation
of $\cp_2$.  When $k=2$ the formulas in Lemma 6 provide an
alternative way of expressing the determinant (\ref{e3}) used
earlier.

\begin{lem}\label{l6}
Let $R= ||f||^2 - ||g||^2 \ge 0$ on the diagonal; then $R \in
\cp_2$ if and only if
\begin{eqnarray}\label{e28}
||(f\wedge g)(z,w)||^2 &\le& ||f(z)||^2 ||f(w)||^2 {\rm sin}^2
\theta_f(z,w) \\
\nonumber &&\qquad \qquad  + ||g(z)||^2 ||g(w)||^2 {\rm sin}^2
\theta_g(z,w)
\end{eqnarray}
for all $z$ and $w$ in $M$. This inequality can also be written as
\begin{eqnarray}\label{e29}
&& ||(f\wedge g)(z,w)||^2 +|\langle f(z),f(w)\rangle|^2 + |\langle
g(z),g(w)\rangle|^2 \\
\nonumber && \qquad \qquad \qquad \qquad  \le ||f(z)||^2
||f(w)||^2 + ||g(z)||^2 ||g(w)||^2.
\end{eqnarray}
\end{lem}

\pf We apply the definition of $\cp_2$; since $R$ is nonnegative
at each point, it is in $\cp_2$ if and only if the determinant
$\Delta_R^2$ is nonnegative at each pair of points $z$ and $w$.
The determinant information yields
\begin{eqnarray}\label{e30}
|R(z, \overline w)|^2 \le R(z, {\overline z}) R(w, {\overline w}).
\end{eqnarray}
Next replace $R(z,{\overline w})$ by its holomorphic
representation in (30). Finally expand, collect terms, and use
(27). It follows that inequality (28) is equivalent to (30).  It
is clear that (28) and (29) are also equivalent. \qed

\medskip

Inequalities (28) and (29) differ from the usual Cauchy-Schwarz
inequality because of the presence of the term $ ||(f\wedge
g)(z,w)||^2$, which relates $f$ and $g$. Inequalities (29) and
(30) are equivalent forms of the nonnegativity of $\Delta_R^2$ as
in (\ref{e3}).

As in the proof of Theorem 1, Lemma 6 is especially useful when
$g$ is scalar-valued. In that case ${\rm sin}^2 \theta_g(z,w)$
vanishes, and formulas (28) and (29) simplify.

\end{document}